\documentclass[11pt]{article}

\usepackage[colorlinks=false,pdfborderstyle={/W 1}]{hyperref}
\def\arXiv#1#2{\href{http://front.math.ucdavis.edu/#1}{{\tt arXiv:#1 [#2]}}} 
\def\arXivo#1{\href{http://front.math.ucdavis.edu/#1}{{\tt [arXiv:#1]}}} 

\usepackage{amsmath}
\usepackage{amssymb}
\usepackage{amsthm}
\usepackage{amsfonts}
\usepackage{graphicx}
\usepackage{mathrsfs}

\input epsf.sty

\newtheorem{thm}{Theorem}[section]
\newtheorem{quest}[thm]{Question}
\newtheorem{conj}[thm]{Conjecture}
\newtheorem{lemma}[thm]{Lemma}
\newtheorem{prop}[thm]{Proposition}

\theoremstyle{remark}\newtheorem{rem}[thm]{Remark}
\theoremstyle{definition}
\newtheorem{defi}{Definition}[section]

\newcommand{\e}{\varepsilon}

\newcommand{\p}{\mathbb{P}}
\newcommand{\E}{\mathbb{E}}
\newcommand{\Z}{\mathbb{Z}}

\newcommand{\Cov}{\mathrm{Cov}}
\newcommand{\Corr}{\mathrm{Corr}}
\newcommand{\Var}{\mathrm{Var}}
\newcommand{\ra}{\rightarrow}
\newcommand{\A}{\mathcal{A}}
\newcommand{\B}{\mathcal{B}}
\newcommand{\C}{\mathcal{C}}

\newcommand{\Maj}{\mathsf{Maj}}
\newcommand{\Piv}{\mathsf{Piv}}
\newcommand{\Inf}{{\bf I}}

\newcommand{\T}{\mathbb{T}}
\newcommand{\lora}{\longrightarrow}

\numberwithin{equation}{section}
\numberwithin{figure}{section}

\let\qqed=\qed
\def\QED{\qqed\medskip}
\let\qed=\QED

\def \eps {\epsilon}

\def \P {{\bf P}}
\def\md{\mid}
\def\Bb#1#2{{\def\md{\bigm| }#1\bigl(#2\bigr)}}
\def\BB#1#2{{\def\md{\Bigm| }#1\Bigl(#2\Bigr)}}
\def\Bs#1#2{{\def\md{\mid}#1(#2)}}

\def\Pso#1{\Bs{\P_{#1}}}
\def\Pbo#1{\Bb{\P_{#1}}}
\def\PBo#1{\BB{\P_{#1}}}
\def\Eso#1{\Bs{\E_{#1}}}
\def\Ebo#1{\Bb{\E_{#1}}}

\def\proofof#1{{ \medbreak \noindent {\bf Proof of #1.} }}

\def\bl{\begin{lemma}}
\def\el{\end{lemma}}
\def\bth{\begin{theorem}}
\def\eth{\end{theorem}}
\def\bc{\begin{corollary}}
\def\ec{\end{corollary}}
\def\bcj{\begin{conjecture}}
\def\ecj{\end{conjecture}}
\def\bpr{\begin{proposition}}
\def\epr{\end{proposition}}
\def\bde{\begin{definition}}
\def\ede{\end{definition}}
\newcommand{\be}{\begin{eqnarray}}
\newcommand{\ee}{\end{eqnarray}}
\newcommand{\bes}{\begin{eqnarray*}}
\newcommand{\ees}{\end{eqnarray*}}

\usepackage{mathrsfs}

\def\1{1\!\! 1}

\title{Noise sensitivity in bootstrap percolation}
\author{Zsolt Bartha and G\'abor Pete}
\date{}

\begin{document}
\maketitle

\begin{abstract} 
Answering questions of Itai Benjamini, we show that the event of complete occupation in 2-neighbour bootstrap percolation on the $d$-dimensional box $[n]^d$, for $d\geq 2$, at its critical initial density $p_c(n)$, is noise sensitive, while in $k$-neighbour bootstrap percolation on the $d$-regular random graph $G_{n,d}$, for $2\leq k\leq d-2$, it is insensitive. Many open problems remain.
\end{abstract}

\section{Introduction and results}

The concept of noise sensitivity of Boolean functions (or events) of i.i.d.~input bits was introduced by Benjamini, Kalai and Schramm in \cite{BKS}, motivated by computer science and statistical mechanics. Roughly speaking, it means that if we resample each input bit with an arbitrarily small probability, then with high probability we cannot predict whether the event occurs in the modified configuration. They proved that a monotone Boolean function is  insensitive to noise if and only if it is positively correlated with a generalized majority function, while their main noise sensitive example was the left-to-right crossing event in critical bond percolation on the square grid. Since then, using different techniques (though all based on discrete Fourier analysis), the exact noise sensitivity of crossing events in planar percolation has been understood well \cite{SS,GPS,AGMT}, and noise sensitivity has become an important area of the analysis of Boolean functions \cite{GS,OD}.  

In the present paper, we study noise sensitivity questions in bootstrap percolation. In this well-known spreading model, one starts with an i.i.d.~Bernoulli$(p)$ set of occupied vertices in a finite or infinite graph $G(V,E)$, then a vertex becomes occupied if at least $k$ of its neighbours are occupied, and this is repeated ad infinitum. The critical density $p_c(G,k)$ is defined as the infimum of initial occupation densities $p$ for which the probability that every vertex becomes occupied is positive or at least $1/2$, for the cases when the graph $G$ is infinite or finite, respectively. 

For infinite transitive graphs, the value of $p_c(G,k)$ is known in a few cases. For $d$-regular trees, $p_c(\T_d,k)$ is the root of an explicit polynomial, it is strictly between $0$ and $1$ for $2\leq k\leq d-1$, and if we take $k_d/d \to\gamma\in[0,1]$, then $\lim_{d\to\infty} p_c(\T_d,k_d)=\gamma$ holds \cite{CLR,BPP}. The recursive computation of $p_c(\T_d,k)$ is based on the fact that incomplete occupation is equivalent to having a vacant $d+1-k$-regular subtree in the initial configuration. For Euclidean lattices, the situation is very different: $p_c(\Z^d,k)=0$ if $k\leq d$, while $p_c(\Z^d,k)=1$ if $k\geq d+1$ \cite{vE,S}. 

For the most natural sequences of finite graphs that converge in the local weak sense \cite{BS,ALy} to the above infinite transitive graphs, the critical densities turn out to converge to their infinite counterparts. For the $d$-regular random graph $G_{n,d}$ on $n$ vertices, Balogh and Pittel \cite{BP} proved that 
\begin{equation}\label{e.Gnd}
\lim_{n\to\infty} p_c(G_{n,d},k) = p_c(\T_d,k)\quad\textrm{ for }2\leq k\leq d-2\,,
\end{equation}
the case $k=d-1$ being somewhat degenerate. For the boxes $[n]^d$, the cases $2\leq k\leq d$ are the interesting ones, and by the work of several people over many years \cite{AL,BPe,CC,CM,H,BBM,BBDCM}, we know that
\begin{equation}\label{e.Znd}
p_c\big([n]^d,k\big) = \left(\frac{\lambda(d,k)+o(1)}{\log_{(k-1)} n}\right)^{d-k+1},\quad\textrm{ as }n\to\infty\,,
\end{equation}
for some constant $\lambda(d,k)>0$, and $\log_{(k-1)}$ denoting $k-1$ times iterated logarithm.

Given the above results, it is natural to wonder about noise sensitivity of bootstrap percolation on a sequence $G_n$ of finite graphs, especially on $[n]^d$ and $G_{n,d}$, as we were explicitly asked by Itai Benjamini. To state our answers precisely,  let us denote by $\C_n$ the event that complete occupation happens in $k$-neighbour bootstrap percolation on $G_n$ with an initial i.i.d.~configuration of density $p_c(n)=p_c(G_n,k)$. This event is determined by the initial configuration $\omega\in\Omega_n=\{-1,1\}^{V(G_n)}$, hence can be identified with a Boolean function $f^\C_n: (\Omega_n,\p_{p_c(n)})\lora\{-1,1\}$, where $1$ stands for being (completely) occupied. We let $\omega^\eps$ denote the configuration obtained from $\omega$ by resampling each bit independently with probability $\eps$, from the measure $\p_{p_c(n)}$. Finally, $\Corr(X,Y)=\Cov(X,Y)/\sqrt{\Var (X)\Var (Y)}$ denotes the correlation between $X$ and $Y$.

\begin{thm}[Random regular graphs]\label{t.Gnd}
For $2\leq k\leq d-2$, the event $\C_n$ on the sequence $G_{n,d}$ is insensitive to noise: for any $\eps>0$ there is a $\delta>0$ such that, for all $n$ large enough, with probability at least $1-\eps$ the random graph $G_{n,d}$ is such that 
$$\Corr_{p_c(n)}(f^\C_n(\omega),f^\C_n(\omega^\eps)) > \delta\,.$$
\end{thm}

\begin{thm}[Euclidean tori and boxes]\label{t.Znd}
For $d\geq 2$ and $k=2$,  the event $\C_n$ on the sequence of tori $\T_n^d=(\Z/n\Z)^d$ and boxes $[n]^d$ is noise sensitive: for any $\eps_n \gg \log\log n / \log n$, we have  
$$\Corr_{p_c(n)}(f^\C_n(\omega),f^\C_n(\omega^{\eps_n})) \to 0\,.$$
\end{thm}

The opposite behaviour in these two cases would call for a general explanation, which we are missing, unfortunately. Recall that we already have a similarly striking difference between the behaviour of the critical densities $p_c(\T_d,k)\in(0,1)$ and $p_c(\Z^d,k)\in\{0,1\}$. In fact, a question formulated in \cite{BPP} is whether a finitely generated group $\Gamma$ is amenable if and only if, for any of its Cayley graphs $G$ and any $k$-neighbour rule, we have $p_c(G,k)\in\{0,1\}$. Such isoperimetric considerations might also play a role in noise sensitivity (see Remark~\ref{r.witness} below), which partly motivates the following general question:

\begin{quest}\label{q.gen}
Is it true that, for any sequence of finite graphs $G_n$, complete occupation by the $k$-neighbour rule at the critical density $p_c(G_n,k)$ is noise sensitive versus insensitive exactly depending on whether $p_c(G_n,k)\to 0$ or $\liminf_n p_c(G_n,k)>0$ holds, respectively? 
\end{quest}

The proof of Theorem~\ref{t.Gnd} is simple, given the work of Balogh and Pittel \cite{BP}: they prove not only~(\ref{e.Gnd}), but also that the threshold window around $p_c(n)$ is as narrow as possible for a monotone event: a width of size $O(1/\sqrt{n})$. This turns out to imply that the event $\C_n$ is correlated with generalized majority at the level $p_c(n)$, and hence is noise insensitive, at least with large $G_{n,d}$-probability. However, this proof leaves the following questions open. For a discussion of the concentration question in part (a), see Remark~\ref{r.concent} below.

\eject
\begin{quest}\label{q.Gnd}
Consider the sequence $G_{n,d}$, and let $2\leq k\leq d-2$.
\begin{itemize}
\item[{\bf (a)}] Is there a $\delta>0$ such that
$$\PBo{G_{n,d}}{\Corr_{p_c(n)}(f_n(\omega),f_n(\omega^\eps)) > \delta}\to 1\,,\qquad\textrm{as }n\to\infty\,?$$
\item[{\bf (b)}] Is the event $\C_n$ actually noise stable? That is, for any $\delta>0$, is there an $\eps>0$ such that, with probability at least $1-\delta$ for all $n$ large enough, or even with probability tending to 1, the random graph $G_{n,d}$ is such that 
$$\Corr_{p_c(n)}(f_n(\omega),f_n(\omega^\eps)) > 1-\delta\,?$$
\end{itemize}
\end{quest}

The result of Theorem~\ref{t.Znd} is motivated by the observation that bootstrap percolation on $[n]^d$ is somewhat similar to one of the simplest examples of a noise sensitive sequence of Boolean functions, \textsf{Tribes}: roughly, we divide the input bits into  logarithmically sized tribes, and the final output is 1 if at least one of these tribes has all its bits equal to 1. Indeed, for $d\geq k=2$, complete occupation in bootstrap percolation is known to be equivalent to the existence of an initially occupied brick-shaped seed of side-lengths $\log^{d-1+o(1)} n$, which then starts growing and occupies the entire box. However, there is also a big difference between bootstrap percolation and \textsf{Tribes}: if we find a $-1$ input bit in a tribe, then we immediately know that this tribe cannot help in getting a 1 at the end, while a sub-box of side-length $\log^{d-1+o(1)} n$, even if many of its bits are $-1$'s, still can be a successful seed. 
As we will discuss in Section~\ref{s.noiseback}, this seems to exclude the so-called ``revealment method'' of Schramm and Steif \cite{SS} to prove noise sensitivity, which does work for \textsf{Tribes}. Instead, we will use the ``squared influence method'' of \cite{BKS}, or more precisely, a generalization of it for input densities tending to 0 \cite{KK,B}. Thus, the following upper bound is one of the key steps in the proof; the lower bound follows immediately from \cite{FK}, as noticed in \cite{BB}:

\begin{prop}\label{prop}
For $d\geq 2$, $k=2$, and $\delta>0$, on the torus $\T_n^d$, let $p(n)$ be such that $\Pso{p(n)}{\C_n}\in (\delta,1-\delta)$. Then, the probability that an initial bit $x$ is pivotal for $\C_n$ (that is, its influence) satisfies 
$$
\frac{\log^{d-o(1)} n}{n^d} \leq \p_{p(n)}(x\text{ is pivotal}) \leq \frac{\log^{d^2-1+o(1)} n}{n^d}\,,
$$
with the $o(1)$ terms depending only on $d$ and $\delta$.
\end{prop}

However, for $d\geq k > 2$, the equivalence between complete occupation and the existence of a nice simple seed is not known (despite \cite{BBM,BBDCM}), hence we do not exactly understand what makes a site pivotal, and hence the following remains open:

\begin{conj}
For any $d\geq k > 2$, the event $\C_n$ on the sequence of tori $\T_n^d$ and boxes $[n]^d$ is noise sensitive.
\end{conj}

Finally, let us remark that a striking application of noise sensitivity in planar percolation is the existence of exceptional times in dynamical percolation on the infinite lattice \cite{SS,GPS}: when the input bits are continuously resampled via independent Poisson clocks, keeping critical percolation as the stationary measure, there exist random times when the configuration has an infinite cluster, even though this event has probability zero at any fixed time. See \cite{St} for a survey of dynamical percolation. Several people have asked independently what happens with bootstrap percolation on $\T_d$ at the critical value $p(\T_d,k)$, when $d+1-k\geq 3$. This question is quite different from exceptional times in planar percolation, because the event of containing a vacant $\ell$-regular subtree in the initial configuration, for $\ell\geq 3$, has a discontinuous phase transition; in this respect, the question is more similar to \cite{PSS}. Based on preliminary work with Marek Biskup, we tend to believe that the answer will be that there do exist exceptional times here, which is quite surprising, given the noise insensitivity on $G_{n,d}$. For the case of $\Z^d$, a consequence of $p_c\in\{0,1\}$ is that there is no meaningful critical infinite system where exceptional times could be studied. However, for bootstrap percolation with more general spreading rules, non-trivial critical densities already occur \cite{BBPS}, hence the question of exceptional times makes sense, and in fact is explicitly asked in that paper. The first question would be: what happens to noise sensitivity there, especially in light of Question~\ref{q.gen}?

\section{Background on noise sensitivity}\label{s.noiseback}

\begin{defi}
A sequence $f_n: (\Omega_n,\p_{p(n)})\lora\{-1,1\}$ of Boolean functions, typically assuming the non-degeneracy property  $\liminf_n \Var_{p(n)}(f_n) > 0$, is called \emph{noise sensitive} if, for every $\eps>0$, 
$$\lim_{n\to\infty} \Corr_{p(n)}(f_n(\omega),f_n(\omega^\eps)) = 0\,.$$ 
The sequence is called \emph{insensitive} if, for every $\eps>0$,
$$\liminf_{n\to\infty} \big|\Corr_{p(n)}(f_n(\omega),f_n(\omega^\eps))\big| = \delta(\eps) > 0\,,$$
and \emph{noise stable} if $\delta(\eps)\to 1$ as $\eps\to 0$.
\end{defi}

Three main ways have been found so far to address noise sensitivity of Boolean functions, all based on discrete Fourier analysis: the first is to apply hypercontractivity, resulting in conditions involving correlation with generalized majority functions and influences \cite{BKS}; the second is the revealment method \cite{SS}; the third is to encode the Fourier expansion as a random subset of the bits, called the Fourier spectral sample, and study the typical size of this random set directly \cite{GPS,OD}. This last approach is possible to do well only in some special cases, hence we will discuss now only the former methods.


\begin{defi}
For any ${\bf w}\in\mathbb{R}^n$ and $s\in\mathbb{R}$, the functions defined on $(\Omega_n,\p_p)$ of the form
\[\Maj_{n,{\bf w},s}(x_1,x_2,\ldots,x_n)=\mathrm{sign}\left(\left(\sum_{i=1}^n w_i\cdot(x_i-(2p-1)) \right)-s\right)\]
are called generalized weighted majority functions.
\end{defi}

Note that $\E_p\left[\sum_{i=1}^n w_i (x_i-(2p-1))\right]=0$.
When $\forall i\in[n]: w_i=1$, we omit ${\bf w}$ from the notation, as well as $s$ if $s=0$. A key theorem is the following:

%
%
%
%
%
%

\begin{thm}[\cite{BKS}]\label{t.maj}
A sequence $\{f_n\}$ of monotone Boolean functions is noise sensitive if and only if
\[\lim_{n\ra\infty}\sup_{{\bf w}\in [0,1]^n}\Corr_p(f_n,\Maj_{n,{\bf w}})=0.\]
\end{thm}

Of course, the easier direction is that generalized majority functions and functions correlated with any of them cannot be noise sensitive, and this is the direction we will use in this paper.

%
%

The second technique we discuss is based on the notion of influences.

\begin{defi}
Given a Boolean function $f:(\Omega_n,\p_{p(n)})\ra\{-1,1\}$, we say that $i\in[n]$ is \emph{pivotal} for $f$ for $\omega$
if $f(\omega)\ne f(\omega^i)$, where for $\omega=(x_1,x_2,\ldots,x_n)$, the configuration $\omega^i$ is defined as $(x_1,\ldots,x_{i-1},-x_i,x_{i+1},\ldots,x_n)$. The \emph{pivotal set} for $f$ is the random set
\[\Piv_f=\{i\in[n]:i\text{ is pivotal for }f\}.\]
\end{defi}

\begin{defi}
For $f$, the \emph{influence} of the $i$th bit is
\[\Inf_i=\Inf_i(f)=\p_{p(n)}(i\in\Piv_f)\,,\]
and the \emph{total influence} of $f$ is
\[\Inf(f)=\sum_{i=1}^n\Inf_i(f).\]
\end{defi}

Another theorem of Benjamini, Kalai and Schramm gives the following criterion for noise sensitivity, in the case $p(n)=1/2$ (other constant densities would also work):

\begin{thm}[\cite{BKS}]\label{inf}
If $f_n: (\Omega_n,\p_{1/2}) \lora\{-1,1\}$, where $\Omega_n=\{-1,1\}^{m_n}$, and
\[\lim_{n\ra\infty}\sum_{i=1}^{m_n}\Inf_i(f_n)^2=0\,,\]
then $\{f_n\}$ is noise sensitive. The converse is also true if $f_n$ is a monotone function for every $n$.
\end{thm}

The influences of some functions can be calculated easily:
\begin{itemize}
\item
For simple majority, $\Maj_n$, we have $\Inf_i=\p\left(\mathrm{Binomial\left(n-1,\frac{1}{2}\right)}=\left\lfloor\frac{n}{2}\right\rfloor\right)\asymp\frac{1}{\sqrt{n}}$ for any $i$, hence Theorem~\ref{inf} implies that it is noise insensitive. In fact, it is easy to see that it is noise stable.
\item 
For $\mathsf{Tribes}_n$ on $n=k2^k$ bits, we have
$\Inf_i=\left(\frac{1}{2}\right)^{k-1}\left(1-\left(\frac{1}{2}\right)^k\right)^{2^k-1}\asymp\frac{1}{2^k}\asymp\frac{\log_2(n)}{n}$
for any $i$. (A bit is pivotal exactly if every other bit in its tribe is $1$ (first factor) and in every other tribe, there is at least
one bit which $-1$ (second factor).) Theorem~\ref{inf} gives noise sensitivity.
\end{itemize}

We will need a generalization of Theorem \ref{inf} that proves noise sensitivity in a case where $p=p(n)$ is allowed to converge to $0$ or to $1$. Two papers, \cite[Theorem 7]{KK}  and  \cite[Theorem 4.1]{B} provide such results, suiting our purposes. We will use the theorem of Keller and Kindler in our calculations, which states the following:

\begin{thm}[\cite{KK}]\label{kk}
With the notation
\[\mathcal{S}_\e(f_n)=\Cov_p(f_n(\omega),f_n(\omega^\eps)),\]
the following holds:
\[\mathcal{S}_\e(f_n)\le(6e+1)\mathcal{W}(f_n)^{\alpha(\e)\cdot\e},\]
where
\[\mathcal{W}(f_n)=p(n)(1-p(n))\sum_{i=1}^{m_n}\Inf_i(f_n)^2,\]
and
\[\alpha(\e)=\frac{1}{\e+\log(2B(p)e)+3\log\log(2B(p)e)}\,,\]
with the so-called hypercontractivity constant
\[B(p)=\frac{\frac{1-p}{p}-\frac{p}{1-p}}{2\log\frac{1-p}{p}}.\]
\end{thm}

Finally, let us say a few words about the revealment method of Schramm and Steif \cite{SS}. We will consider algorithms that ask the values of some of the input bits one after the other, and stop when, from the values of the bits revealed, the value of $f$ can be determined. During this computation, the algorithm can decide which bit to query next based on the previously asked values, and may also use auxiliary randomness. For such a {randomized algorithm} $A$, let $Q_A$ denote the random set of queried bits; it depends on the randomness of the actual value of the bits as well as on the auxiliary randomness used during the running of the algorithm.

\begin{defi}
The \emph{revealment of a randomized algorithm} $A$ for a Boolean function $f$ is
\[\delta_A:=\max_{i\in[n]}\p(i\in Q_A).\]
The \emph{revealment of a Boolean function} $f$ is
$$\delta_f:=\inf_A\delta_A,$$
where the infimum is taken over all possible randomized algorithms that compute $f$.
\end{defi}

The following theorem shows how revealment enters the study of noise sensitivity. It may seem completely counterintuitive: if the function can be computed by revealing only few bits, then it is noise sensitive, i.e., resampling few bits will destroy all information.

\begin{thm}[\cite{SS}]\label{revthm}
If, for a sequence of Boolean functions $\{f_n\}$, the revealments satisfy
\[\lim_{n\ra\infty}\delta_{f_n}=0,\]
then $\{f_n\}$ is noise sensitive. (Their result also gives quantitative information, but let us give just this weak version now.)
\end{thm}

As two examples, it is easy to see that majority cannot be computed with small revealment, while for $\mathsf{Tribes}$, the following natural algorithm has small revealment: we query the bits of a tribe one-by-one until we find all of them to be $1$, and the output is 1, or we find a $-1$ in the tribe, at which point we move to the next tribe.

It would be nice if the converse to Theorem~\ref{revthm} was true, and noise sensitive functions always had algorithms with small revealment. However, this is false, as pointed out in \cite[Section VIII.6]{GS}, by the example of clique containment in the Erd\H{o}s-R\'enyi random graph $G(n,1/2)$. We suspect that $2$-neighbour bootstrap percolation on $[n]^2$ is very similar to this example, having uniformly positive revealment, even though it is noise sensitive by our Theorem~\ref{t.Znd}. Having no converse is also a pity for the following reason:

\begin{rem}\label{r.witness}
Consider $k$-neighbour bootstrap percolation on a finite $d$-regular expander graph $G(V,E)$, i.e., a graph with \emph{edge Cheeger constant} 
$$h_E:=\min \left\{\frac{|\partial_E S|}{|S|} : S\subset V,\  |S|\leq |V|/2 \right\} > 0\,;$$
think, e.g., of $G_{n,d}$, which is known to be an expander with high probability; see, e.g., \cite[Section 1.2]{Lub}. Assume that $d<2k+h_E$, just like in \cite[Theorem 1.4]{BPP}, which concerned bootstrap percolation on nonamenable infinite graphs. 
The edge boundary of the occupied set can always increase by at most $d-2k$ during the $k$-neighbour occupation process, as we occupy vertices one-by-one. Thus, if an initial set $A$ of size less than $|V|/2$ occupies all of $G$, then during this process there is a moment where the number of occupied vertices is exactly $|V|/2$ (assuming that $|V|$ is even, for notational simplicity), at which point the edge boundary is at least $h_E  |V|/2$, hence the initial edge boundary $|\partial_E A|$, which is clearly at most $d|A|$ on the one hand, was of at least size $h_E |V|/2 - ( |V|/2-|A|)(d-2k)$ on the other hand. This is possible only if $|A| \geq (h_E+2k-d)|V|/(4k)$. Therefore, any subset of the input bits that witnesses successful occupation must have size at least $c|V|$, with $c=(h_E+2k-d)/(4k)$, which certainly implies that no algorithm with revealment smaller than $c$ can exist in this case. By the above examples, this does not automatically imply that complete occupation is noise insensitive in this case, but we do conjecture that the condition $d<2k+h_E$ does imply that. In particular, this would give a more robust reason, independent of \cite{BP}, why bootstrap percolation on $G_{n,d}$, at least for $k>d/2$, is noise insensitive.
\end{rem}

\section{The case of $G_{n,d}$}

We start by recalling what Balogh and Pittel \cite{BP} proved about the critical probability and transition window for $G_{n,d}$,
a graph chosen uniformly at random from the $d$-regular graphs ($d\ge 3$) on $n$ vertices. (It is proven \cite{W1} that $G_{n,d}$ is connected with probability tending to $1$ as $n\rightarrow\infty$.) We set $A$ to be the set of initially occupied vertices. At this point, let the density of initially occupied vertices be fixed: $p=|A|/n$. We can choose the vertices in $A$ deterministically, and then decide about the edges of the graph to obtain a uniformly random $d$-regular graph with $pn$ initially occupied vertices. Denote by $I_f$ the set of eventually vacant vertices, following $k$-neighbour bootstrap percolation, with $k\ge 2$ fixed.

Introduce
\begin{equation}\label{p*}
p^*:=1-\inf_{y\in(0,1)}R(y),\quad R(y):=\frac{y}{\p(\mathrm{Binomial}(d-1,1-y)<k)}.
\end{equation}
For $k<d-1$ the infimum is attained at an interior point $y^*$.
For example if $k=2$ and $d>3$, then
\[y^*=\frac{(d-1)(d-3)}{(d-2)^2},\quad p^*=1-\frac{(d-2)^{2d-5}}{(d-1)^{d-2}(d-3)^{d-3}},\]
while for $k=d-1$
\[p^*=1-\frac{1}{(d-1)}.\]

By the definition of $p^*$, the equation
\begin{equation}\label{py}
p=1-R(y)
\end{equation}
has no root $y\in(0,1)$ if $p>p^*$.
It turns out that for $p<p^*$ \eqref{py} has exactly one root $\hat{y}(p)$ in $(0,1)$ if $k=d-1$, and exactly
two roots in $(0,1)$ if $k<d-1$, in which case it is the larger root that we denote by $\hat{y}(p)$.

With these notations, the main theorem of Balogh and Pittel is the following.

\begin{thm}\label{BP1}
Suppose that $2\le k<d-1$. Let $\omega=\omega(n)\rightarrow\infty$ so slowly that $\log\omega(n)=o(\log n)$.

(i) If $p>p^*+\omega n^{-1/2}$, then
\[\lim_{n\rightarrow\infty}\p(I_f=\emptyset)=1.\]

(ii) If $p<p^*-\omega n^{-1/2}$, then for every $\gamma>1/2$
\[\lim_{n\rightarrow\infty}\p\left(|I_f|=nh(\hat{y}(p))+O(n^{1/2}\omega^{\gamma}(p^*-p)^{-1/2})\right)=1,\]
\quad\quad\quad where
\[h(y):=(1-p)\p(\mathrm{Binomial}(d,1-y)<k).\]
\end{thm}

Thus, for $k<d-1$, $p^*$ is the critical probability and the transition window has size at most $\omega n^{-1/2}$ with a sequence
$\omega(n)$ tending to $\infty$ arbitrarily slowly, so it has size of order $n^{-1/2}$ at most. Part (ii) means intuitively that
for $p<p^*$ outside the transition window, the numbers of the eventually occupied neighbours of an eventually vacant vertex
is distributed binomially in the limit with success probability $1-\hat{y}(p)$.

The case $k=d-1$ is quite different. Now any cycle in $G_{n,d}$ is a fort: if it is vacant in the initial configuration, then it remains vacant forever. And, as proven in \cite{W2}, $G_{n,p}$ does contain cycles
(of length $3$) with positive limiting probability. It follows that in the initial configuration we have a completely vacant fort
with positive limiting probability if $p<1$, therefore $\limsup_{n\rightarrow\infty}\p($complete occupation in $G_{n,d})<1$ for every such $p$.
Now the theorem proven by Balogh and Pittel is the following.

\begin{thm}\label{BP2}
(i) Suppose that $p\ge p^*+n^{-\e}$, where $\e=\e(n)\searrow 0$ and $\e\log(n)\rightarrow\infty$. Then
\[\lim_{n\rightarrow\infty}\p(|I_f|=O((p-p^*)^{-3/2}))=1.\]
(ii) Let $p\le p^*-\omega n^{-\sigma}$, $\sigma=\frac{1}{2d +5}$, $\omega=\omega(n)\rightarrow\infty$, $\log\omega=o(\log n)$. Then
\[\lim_{n\rightarrow\infty}\p(|I_f|=nh(\hat{y}(p))+O(n^{1-3\sigma}\omega^{6\sigma}(p^*-p)^{-1}))=1.\]
(The remainder term is negligible compared to $nh(\hat{y}(p))$.)
\end{thm} 

Thus, $p^*$ is the critical probability for the $k=d-1$ case too, but the transition window seems asymmetric and much wider than for $k<d-1$.

So far, the size of the set $A$ (or equivalently, the set $A$ itself) was deterministic. But we can easily extend the above theorems for
a class of random initially occupied vertex sets:

\begin{thm}\label{BP3}
Suppose that the set $A$ of initially occupied vertices is random, independently of $G_{n,d}$, satisfying the following condition:
\[\lim_{n\rightarrow\infty}\p(n^{-1}\left||A|-\E|A|\right|\le\lambda n^{-1/2})=1\]
for any $\lambda=\lambda(n)\rightarrow\infty$, i.e., $|A|=\E|A|+O_p(n^{1/2})$. Then the statements of Theorems \ref{BP1} and \ref{BP2} hold
with $p:=n^{-1}\E|A|$.
\end{thm}

What is important for us is that the condition for $|A|$ in Theorem \ref{BP3} holds
when each vertex is occupied independently with probability $p$.

%
%

\proofof{Theorem~\ref{t.Gnd}}
As mentioned before, we will show that $f^\C_n$ is uniformly correlated with a certain generalized majority function because of the tiny transition window, and then we can use Theorem~\ref{t.maj}.

We denote the probability measure on the regular graphs by $\p_{G_{n,d}}$, and the product measure of this with $\p_{p(n)}$ by $\p^G_{p(n)}$. Firstly, we deduce from Theorem \ref{BP3} that for any $\e\in(0,\frac{1}{2})$, there exists some constant $K_1>0$ such that for any $n$
\begin{equation}\label{K1}
\p^{G}_{p(n)}(\C_n)\ge 1-\e\quad\text{whenever}\quad p(n)>p^*+K_1n^{-1/2}.
\end{equation}
To this end, assume that $p_{1-\e}(n)$ (the probability for which
$\p^{G}_{p_{1-\e}(n)}(\C_n)=1-\e$) is equal to $p^*+K(n)n^{-1/2}$, where $K(n)\ra\infty$. Then $\sqrt{K(n)}$ tends to
infinity as well, and for $q(n)=p^*+\sqrt{K(n)}n^{-1/2}$, $\p^{G}_{q(n)}(\C_n)\le 1-\e$ for every $n$, thus, it does not tend to $1$, which
contradicts Theorem \ref{BP3}. (We can similarly rule out  the case that infinity is just an accumulation point of $\{K(n)\}$ by using the appropriate subsequence $\{n_\ell\}$ in the proof.) So, $K_1=\sup_n K(n)$ satisfies \eqref{K1}.

Proving a lower bound similarly, we get that there exists a constant $K_2>0$ such that, for any $n$,
\begin{equation}\label{K2}
|p_c(n)-p^*|\le K_2 n^{-1/2}.
\end{equation}
Consequently, Theorem \ref{BP1} remains true if we substitute $p_c$ in the place of $p^*$.
Using this, and the same reasoning as in the proof of \eqref{K1} above, we get that for any $\e\in(0,\frac{1}{2})$ there is a constant $K>0$ such that, for any $n$,
\begin{equation}\label{K}
\p^{G}_{p(n)}\big(\C_n\,|\,\#A_n=m\big)\ge 1-\e\quad\text{whenever}\quad m>np_c(n)+Kn^{1/2}=:r(K,n),\,m\in\Z,
\end{equation}
where $\#A_n$ is the number of initially occupied vertices. (Note that the probability on the left-hand side does not depend on $p(n)$.)

With this $K$, from \eqref{K} we have that for any $n$
\begin{align}\label{K'}
\p^{G}_{p_c(n)}&\big(\C_n\,|\,\#A_n>r(K,n)\big)\notag\\
&=\frac{\p^{G}_{p_c(n)}\big(\C_n\,,\,\#A_n>r(K,n)\big)}{\p^{G}_{p_c(n)}(\#A_n>r(K,n))}\notag\\
&=\frac{\sum_{m>r(K,n)}\p^{G}_{p_c(n)}\big(\C_n\,|\,\#A_n=m\big)\p^{G}_{p_c(n)}(\#\{A\}=m)}
{\sum_{m>r(K,n)}\p^{G}_{p_c(n)}(\#\{A\}=m)}
\notag\\
&\ge\frac{\sum_{m>r(K,n)}(1-\e)\p^{G}_{p_c(n)}(\#\{A\}=m)}{\sum_{m>r(K,n)}\p^{G}_{p_c(n)}(\#\{A\}=m)}\notag\\
&=1-\e.
\end{align}

Since for any $\e\in(0,\frac{1}{2})$ there is a $K$ such that \eqref{K'} holds, the following is true: for any $\e\in(0,\frac{1}{2})$ there
exist a constant $L>0$ such that for any $n$
\begin{equation}\label{L}
\p_{G_{n,d}}\big(\p_{p_c(n)}\big(\C_n\,|\,\#A_n>r(L,n)\big)>1-\e\big)>1-\e.
\end{equation}
Using the notation $A_n=(x_1^{(n)},x_2^{(n)},\ldots,x_n^{(n)})\in\Omega_n$, the condition $\{\#A_n>r(L,n)\}$
is
\[\frac{\sum_{i=1}^n x_i^{(n)}+n}{2}>p_c(n)n+Ln^{1/2},\]
or
\[\sum_{i=1}^n (x_i^{(n)}-2p_c(n)+1)>2Ln^{1/2}.\]
Therefore \eqref{L} means exactly that, with $\p_{G_{n,d}}$-probability greater than $1-\e$, $f^\C_n$ is
correlated with the generalized majority function $\Maj_{n,2Ln^{1/2}}$ on $\Omega_n$ with probability parameter $p(n)$,
and this correlation is uniformly larger than $0$. This is because,
by the properties of the binomial distribution,
\[\lim_{n\ra\infty}p_{2L}(n)=\lim_{n\ra\infty}\p_{p_c(n)}(\Maj_{n,2Ln^{1/2}}=1)=p_{2L}>0,\]
and
\begin{align*}
\p_{p_c(n)}\big(\C_n\cap&\#A_n>r(L,n)\big)=\\
&=\p_{p_c(n)}\big(\C_n\,|\,\#A_n>r(L,n)\big)\p_{p_c(n)}\big(\#A_n>r(L,n)\big)\\
&\ge(1-\e)\p_{p_c(n)}\big(\#A_n>r(L,n)\big)\\
&=\p_{p_c(n)}(\C_n)\p_{p_c(n)}\big(\#A_n>r(L,n)\big)
+\left(\frac{1}{2}-\e\right)\p_{p_c(n)}\big(\#A_n>r(L,n)\big),
\end{align*}
so
$$
\Ebo{p_c(n)}{f^\C_n\, \Maj_{n,2Ln^{1/2}}}-\Ebo{p_c(n)}{f^\C_n} \, \Ebo{p_c(n)}{\Maj_{n,2Ln^{1/2}}}\ge
\left(\frac{1}{2}-\e\right)p_{2L}(n)\,,
$$
with $\p_{G_{n,d}}$-probability greater than $1-\e$ for any $n$. That is, for $n$ large enough, we have
$$
\Cov_{p_c(n)}\big(f^\C_n\,,\,\Maj_{n,2Ln^{1/2}}\big)\ge\left(\frac{1}{2}-\e\right)\frac{p_{2L}}{2}>0
$$
with probability at least $1-\eps$, and Theorem~\ref{t.maj} proves the noise instability of $\{f^\C_n\}$.
\qed


\begin{rem}\label{r.concent}
Consider the function $\varphi(G):=\p_{p_c(n)}\big(\C_n\,|\,\#A_n>r(K,n)\big)$, for bootstrap percolation on a given graph $G$. We saw above that $\Ebo{G_{n,d}}{\varphi(G)} > 1-\eps$ if $K$ is large enough. To answer Question~\ref{q.Gnd}~(a), one would need to prove that $\varphi(G)$ is concentrated around its mean: for instance, that $\Pbo{G_{n,d}}{\varphi(G) > 1-2\eps}\to 1$ as $n\to\infty$. A natural idea would be to try and use Azuma-Hoeffding for Lipschitz martingale sequences (see, e.g., \cite[Section 1.2]{PGG}), where $G_{n,d}$ is explored edge-by-edge, or something similar in the configuration model. For this to give the concentration we want, we would need that the sum of the squares of the step-by-step Lipschitz constants is $o(1)$, so we would basically need that changing one edge in $G\sim G_{n,d}$ changes $\varphi(G)$ typically by $o(1/\sqrt{n})$. Being in a bounded degree graph, this seems to be very similar to proving that the probability that a vertex is pivotal is $o(1/\sqrt{n})$. However, we know that this is false, since by Theorem~\ref{inf} it would give noise sensitivity for complete occupation.
\end{rem}


\section{The case of Euclidean boxes and tori
}

The following two simple lemmas will be of key importance for us. A subset $S$ of the vertices is called \emph{internally spanned} if it gets occupied by the process restricted to $S$, without the help of occupied vertices outside $S$.

\begin{lemma}[\cite{AL}]\label{AL2}
In a successful bootstrap percolation process with parameter $k=2$ on $[n]^d$, $d\ge 2$, for every integer $\ell\in[1,n]$ there is an integer $m\in[\ell,2\ell]$ such that $[n]^d$ contains an internally spanned rectangle with longest side $m$. The same holds for $\T_n^d$.
\end{lemma}

The second lemma is a useful comparison between the process on the torus $\T_n^d$ and on $[n]^d$. Our arguments below will work better for the transitive case $\T_n^d$, and then we will translate them to the case of $[n]^d$ using this lemma.

\begin{lemma}[\cite{BB}]\label{lemBB}
For any $d\geq 2$, $k=2$, and $p\in(0,1)$,
$$\Pbo{p}{\C_{[n]^d}}<\Pbo{p}{\C_{\T_n^d}}<\Pbo{p+n^{-1/3}}{\C_{[n]^d}}+2n^{-1/3}\,.
$$
\end{lemma}

Let us stress that this lemma is nothing mysterious; in fact, one could prove it using ideas similar to the ones below.


\proofof{Proposition~\ref{prop}} We will present the proof in detail for $d=2$, then briefly explain what needs to be changed for $d>2$.

Let $\lambda:=\Pso{p(n)}{\C_{\T_n^2}}\in (\delta,1-\delta)$, and fix a vertex $x\in \T_n^2$. We want to give an upper bound for
\[{\bf I}_x=\p_{p(n)}(x\text{ is pivotal for the event of complete occupation}).\]
Assume that there is complete occupation. Then, by Lemma \ref{AL2}, there is an internally spanned rectangle in $\T_n^2$ with longest side
$\ell\in[f(n),2f(n)]$, where $f(n)=\lfloor c_1 \log^{1+\e} n\rfloor$ with a fixed $\e>0$ and constant $c_1$ specified later. If there are more such rectangles,
we choose one of them uniformly at random (independently of everything else); let this be $R$. The transitivity of $\T_n^2$ implies that
$$
\Pbo{p(n)}{x\notin R \md \C_n} \ge 1-\frac{4f(n)^2}{n^2},
$$
since at most $(2f(n))^2$ vertices can be the bottom left corner of $R$ if it contains $x$.
Thus, if we denote by $\A_x$ the event that there is an internally spanned rectangle in $\T_n^2$ with longest side
$\ell\in[f(n),2f(n)]$ that does not contain $x$, we have
$$
\p_{p(n)}(\A_x) \geq \p_{p(n)}(\A_x\cap\C_n)\ge \p_{p(n)}(\C_n)\left(1-\frac{4f(n)^2}{n^2}\right)=\lambda\left(1-\frac{4f(n)^2}{n^2}\right).
$$

On the event $\A_x$, we choose a rectangle $R$ as above, uniformly at random if there are more than one. On this event, let us consider another event, $\A'_R$, that this internally spanned $R$ continues to grow without using the possibly initially occupied $x$ until it becomes a square with side length $g(n)=\lfloor c_2 \log^2 n\rfloor$ with a constant $c_2$ specified later. This $\A'_R$ clearly happens if any $f(n)$-segment
(i.e., rectangle of size $1\times f(n)$ or $f(n)\times 1$) with distance at most $g(n)$ from $R$ and not containing $x$ has at least one initially occupied vertex in it. Since there are at most $Cg(n)^2$ such $f(n)$-segments (with some fixed constant $C$), this event has probability
$$
\Pbo{p(n)}{\A'_R \md \A_x} \ge 1-Cg(n)^2(1-p(n))^{f(n)}.
$$
Thus, if we denote by $\A'_x$ the event that $\A_x$ and $\A'_R$ both occur, then 
\[\p_{p(n)}(\A'_x)\ge \lambda\left(1-\frac{4f(n)^2}{n^2}\right) \left(1-Cg(n)^2(1-p(n))^{f(n)}\right) 
\ge \lambda - \frac{4f(n)^2}{n^2}-
Ce^{2\log g(n)-p(n)f(n)}.\]

Finally, consider the event $\B_x$ that all $g(n)$-segments that do not
contain $x$ have at least one initially occupied vertex. The number of such $g(n)$-segments is at most $2n^2$, therefore a union bound gives
\[\p_{p(n)}(\B_x)\ge 1-2n^2(1-p(n))^{g(n)}\ge 1-2e^{2\log n-p(n)g(n)}.\]

It is easy to see that if both $\A'_x$ and $\B_x$ occur, then we get complete occupation without using that $x$ might be initially occupied. Another union bound gives
\begin{align}\label{eqab}
\p_{p(n)}(\A'_x\cap\B_x)&\ge\lambda-\frac{4f(n)^2}{n^2}-Ce^{2\log g(n)-p(n)f(n)}-2e^{2\log n-p(n)g(n)}\notag\\
&=\lambda-\frac{4\lfloor c_1\log^{1+\e}n\rfloor^2}{n^2}-Ce^{2\log (\lfloor c_2 \log^2 n\rfloor)-p(n)\lfloor c_1 \log^{1+\e} n\rfloor}-2e^{2\log n-p(n)\lfloor c_2 \log^2 n\rfloor}.
\end{align}

Here we can see how we have to choose the constants in order for $\lambda-\frac{4\lfloor c_1 \log^{1+\e}n\rfloor^2}{n^2}$ to be the dominant term. We know from (\ref{e.Znd}) that $p(n)=\Theta\left(\frac{1}{\log n}\right)$. Therefore, we can choose $c_1$ and $c_2$ large enough for the last two terms in \eqref{eqab} to be negligible compared to $\frac{2 \lfloor c_1 \log^{1+\e}n\rfloor^2}{n^2}$. Thus,
\begin{align*}
\p_{p(n)}(\text{complete occupation regardless of the initial status of }x)
=\lambda-O\left(\frac{\log^{2(1+\e)} n}{n^2}\right),
\end{align*}
as $n\ra\infty$. This is true for any $\e>0$, hence
\begin{equation}\label{piv1}
\p_{p(n)}(\C\cap\{x\text{ is pivotal}\})=O\left(\frac{\log^{2+o(1)} n}{n^2}\right).
\end{equation}

Consider the mapping
\[\mathcal{M}_x:\Omega_{n^2}\ra\Omega_{n^2},\ \omega\mapsto\omega^x,\]
which flips the initial state of $x$. This is clearly a bijection, and
\[\mathcal{M}_x(\C\cap\{x\text{ is pivotal}\})=\neg\C\cap\{x\text{ is pivotal}\}.\]
In $\C\cap\{x\text{ is pivotal}\}$ $x$ is initially occupied, and in $\neg\C\cap\{x\text{ is pivotal}\}$ $x$ is initially vacant.
Thus,
\begin{equation}\label{piv2}
\p_{p(n)}(\C\cap\{x\text{ is pivotal}\})=\Theta\left(\frac{1}{\log n}\right)\p_{p(n)}(\neg\C\cap\{x\text{ is pivotal}\}).
\end{equation}

From \eqref{piv1} and \eqref{piv2} we get
\[{\bf I}_x=\p_{p(n)}(x\text{ is pivotal})=O\left(\frac{\log^{3+o(1)} n}{n^2}\right)\,,\]
finishing the proof of the proposition for $d=2$.

For $d>2$, the argument is the same, except that we have to take $f(n)=\lfloor c_1 \log^{d-1+\e}n\rfloor$ and $g(n)=\lfloor c_2 \log^d n\rfloor$, and have to use that $p(n)=\Theta\left(\frac{1}{\log^{d-1}n}\right)$.
\qed
\medskip

\proofof{Theorem~\ref{t.Znd}} Again, let us focus on $d=2$. Since Proposition~\ref{prop} holds for all $x\in\T_n^2$, we have
\begin{equation}\label{suminf}
\sum_{x\in\T_n^2}\Inf_x^2=n^2 O\left(\frac{\log^{6+o(1)} n}{n^4}\right)=O\left(\frac{\log^{6+o(1)} n}{n^2}\right).
\end{equation}

Now we use Theorem \ref{kk}, for $p=p_c(n)$, to show that noise sensitivity follows:
\[B(p)=\frac{\frac{1-p}{p}-\frac{p}{1-p}}{2\log\frac{1-p}{p}}=\Theta\left(\frac{\log n}{\log\log n}\right),\]
from which
\[\alpha(\e)=\frac{1}{\e+\log(2B(p)\e)+3\log\log(2B(p)e)}=\Theta\left(\frac{1}{\log\log n}\right).\]
From \eqref{suminf} we have
\[\mathcal{W}(\C_{\T_n^2})=p(1-p)\sum_{x\in\T_n^2}\Inf_x^2=O\left(\frac{\log^{5+o(1)} n}{n^2}\right).\]
So, for the noise stability we get the following bound:
\[\mathcal{S}_\e (\C_{\T_n^2})\le(6e+1)\mathcal{W}(\C_{\T_n^2})^{\alpha(\e)\cdot\e}
=O\left(\left(\frac{\log^{5+o(1)} n}{n^2}\right)^{\e\cdot\frac{1}{\log\log n}}\right)\]
To see that the right-hand side tends to $0$ when the noise is $\eps=\eps_n \gg \log\log n/\log n$, we prove that its logarithm tends to $-\infty$:
\[\log\left[\left(\frac{\log^{5+o(1)} n}{n^2}\right)^{\e\cdot\frac{1}{\log\log n}}\right]=
\frac{\e}{\log\log n}\left[(5+o(1))\log\log n-2\log n\right]\to -\infty,\]
which shows noise sensitivity for the case of the torus $\T_n^2$.

We now want to deduce from this the analogous result for $[n]^2$. As a preparation, still working in the torus, note that Russo's formula (see, e.g., \cite[Section 12.3]{PGG}), combined with Proposition~\ref{prop}, gives
\begin{equation}\label{Russo}
\log^{2-o(1)} n \leq \frac{d}{dp}\Pso{p}{\C_{\T_n^2}} = \Ebo{p}{ |\Piv(\C_{\T_n^2})| } \leq \log^{3+o(1)} n\,,
\end{equation}
for any $p=p(n)$ that satisfies $\Pso{p}{\C_{\T_n^2}} \in (\eps,1-\eps)$. In particular, for such sequences $p(n)$, we have
$$
\Pso{p(n)}{\C_{\T_n^2}} < \Pso{p(n)+n^{-1/3}}{\C_{\T_n^2}} < \Pso{p(n)}{\C_{\T_n^2}} + \frac{\log^{3+o(1)} n}{n^{1/3}}\,.
$$
By Lemma~\ref{lemBB}, this implies that 
\begin{equation}\label{probdiff}
\Pbo{p(n)}{\C_{\T_n^2}} - \frac{\log^{3+o(1)} n}{n^{1/3}} < \Pbo{p(n)}{\C_{[n]^2}} < \Pbo{p(n)}{\C_{\T_n^2}}\,.
\end{equation}
From this, using now the lower bound in (\ref{Russo}), we also get
\begin{equation}\label{pcdiff}
p_c(\T_n^2) < p_c([n]^2) <  p_c(\T_n^2) + \frac{\log^{1+o(1)}}{n^{1/3}}\,.
\end{equation}

Now, we have proved above that 
\begin{equation}\label{Tcov}
\lim_{n\ra\infty}\Cov_{p_c(\T_n^2)} \left( f^\C_{\T_n^2}(\omega),N_{p_c(\T_n^2)}^\eps f^\C_{\T_n^2}(\omega) \right)=0\,,
\end{equation}
for any $\eps=\eps_n \gg \log\log n / \log n$, where $N_p^\eps f(\omega):=\Eso{p}{f(\omega^\eps) \md \omega}$ is the noise operator. By \eqref{probdiff} and \eqref{pcdiff}, we have that $f^\C_{\T_n^2}$ is close to $f^\C_{[n]^2}$ both under the product measure with density $p_c(\T_n^2)$ and with density $p_c([n]^2)$. Furthermore, since $N_p^\eps$ is a contraction in $L^2(\Omega_n,\P_p)$ for any $p$, the same holds for $N^\eps f^\C_{\T_n^2}$ and $N^\eps f^\C_{[n]^2}$. Therefore, (\ref{Tcov}) implies 
$$
\lim_{n\ra\infty}\Cov_{p_c([n]^2)} \left( f^\C_{[n]^2}(\omega),N_{p_c([n]^2)}^\eps f^\C_{[n]^2}(\omega) \right) =0\,,
$$
which completes our proof for $[n]^2$.

For the case of $\T_n^d$ and $\Z_n^d$, we just have to use the $d$-dimensional versions of Proposition~\ref{prop} and Lemma~\ref{lemBB}.
 \qed

\section{Acknowledgments}

We are indebted to Itai Benjamini for asking us these questions, and to Rob Morris for suggesting a polylog improvement to our original proof in the case of Euclidean lattices.

This research was partially supported by the Hungarian National Science Fund OTKA grant K109684, and by the MTA R\'enyi Institute ``Lend\"ulet'' Limits of Structures Research Group.

\vskip 0.8 cm

\noindent {\bf Zsolt Bartha}\\
Department of Statistics, University of California, Berkeley\\
{\tt bartha `at' berkeley `dot' edu}

\vskip 0.4 cm

\noindent {\bf G\'abor Pete}\\
R\'enyi Institute, Hungarian Academy of Sciences, Budapest,\\
and Institute of Mathematics, Budapest University of Technology and Economics\\
\url{http://www.math.bme.hu/~gabor}\\
{\tt gabor `at' math `dot' bme `dot' hu}

\end{document}